\documentclass[11pt]{amsart}

\usepackage{amsmath,amssymb,amsthm,amsfonts,amscd,flafter,epsf, epsfig,graphicx}
\usepackage[all]{xy}

\title{Comultiplicativity of the Ozsv{\'a}th-Szab{\'o} contact invariant}
\author[John A. Baldwin]{John A. Baldwin}
\address {Department of Mathematics, Columbia University\\ New York, NY 10027}
\email {baldwin@math.columbia.edu}
\date{}

\newcommand\cheeg{\widehat{CF}}
\newcommand\heeg{\widehat{HF}}

\newcommand\ZZ{\mathbb{Z}_2}
\newcommand\Sig{S}
\newcommand\xg{\overline{x_g}}
\newcommand\xh{\overline{x_h}}
\newcommand\xhg{\overline{x_{hg}}}
\newcommand\xgi{(\overline{x_g})_i}
\newcommand\xhi{(\overline{x_h})_i}
\newcommand\xhgi{(\overline{x_{hg}})_i}
\newcommand\cf{\widehat{CF}}
\newcommand\hf{\widehat{HF}}
\newcommand\cfl{CF^{\leq 0}}

\newcommand\cfp{CF^+}
\newcommand\hfp{HF^+}
\newcommand\Xgh{X_{g,h}}
\newcommand\Yh{Y_h}
\newcommand\Yg{Y_g}
\newcommand\Yhg{Y_{hg}}
\newcommand\x{\otimes}
\newcommand\s{\mathfrak{s}}
\newcommand\sg{\mathfrak{s}_g}
\newcommand\sh{\mathfrak{s}_h}
\newcommand\shg{\mathfrak{s}_{hg}}
\newcommand\al{\bar{\alpha}}
\newcommand\be{\bar{\beta}}
\newcommand\ga{\bar{\gamma}}

\newtheorem{thm}{Theorem}[section]

\newtheorem{claim}[thm]{Claim}

\newtheorem{definition}[thm]{Definition}

\newtheorem{cor}[thm]{Corollary}

\begin{document}
\maketitle
\begin{abstract}  
Suppose that $\Sig$ is a surface with boundary and that $g$ and $h$ are diffeomorphisms of $\Sig$ which restrict to the identity on the boundary. Let $Y_g,$ $Y_h$, and $Y_{hg}$ be the three-manifolds with open book decompositions given by $(S,g)$, $(S,h)$, and $(S,hg)$, respectively. We show that the Ozsv{\'a}th-Szab{\'o} contact invariant is natural under a comultiplication map $\tilde{\mu}:\heeg(-Y_{hg}) \rightarrow \heeg(-Y_{g}) \otimes \heeg(-Y_{h}). $ It follows that if the contact invariants associated to the open books $(\Sig, g)$ and $(\Sig, h)$ are non-zero then the contact invariant associated to the open book $(\Sig, hg)$ is also non-zero. We extend this comultiplication to a map on $\hfp(-\Yhg)$, and as a result we obtain obstructions to the three-manifold $\Yhg$ being an $L$-space. We also use this to find restrictions on contact structures which are compatible with planar open books.
\end{abstract} 

\section{Introduction}
In 2002, Ozsv{\'a}th and Szab{\'o} discovered an invariant of contact 3-manifolds constructed as follows \cite{osz1}. Given a contact 3-manifold $(Y,\xi)$, we can find a compatible (in the sense of Giroux \cite{giroux}) open book decomposition $(\Sig, g)$ of $Y$ with connected binding $K \subset Y$. If the genus of $\Sig$ is $k$ then the knot Floer homology, $\widehat{HFK}(-Y,K,-k) \cong \mathbb{Z}$. Moreover, it is possible to find a Heegaard diagram for which the knot Floer homology in this $-k$ th filtration level is generated by an element of the chain complex $\widehat{CFK}(-Y,K)$ which also represents a cycle in the chain complex $\widehat{CF}(-Y)$. The contact invariant $c(\xi)$ is defined to be the image of this cycle in $\widehat{HF}(-Y)$. 

The class $c(\xi)$ is well-defined up to sign (we use $\ZZ$ coefficients throughout to avoid ambiguity in sign), and it is an invariant of $ \xi$ up to isotopy of the contact structure. This invariant encodes information related to the tightness of $\xi$. For instance, Ozsv{\'a}th and Szab{\'o} prove that if $\xi$ is overtwisted, then $c(\xi) = 0$. On the other hand, if $\xi$ is Stein fillable or strongly fillable, then $c(\xi) \neq 0$ \cite{osz1}, \cite{osz8}. In a previous paper we show, in the case of contact structures compatible with genus one, one boundary component open books, that $c(\xi) = 0$ if and only if $\xi$ is overtwisted for all but a small family of open books with reducible monodromies \cite{bald1}. Honda, Kazez, and Mati{\'c} have since shown this to be the case for all such open books \cite{hkm3}, \cite{hkm1}, \cite{hkm2}. However, the precise relationship between $c(\xi)$ and the tightness of $\xi$ is still unknown -- there are tight contact structures with vanishing contact invariant \cite{ghiggini}. In fact, Lisca and Stipsicz conjecture that the contact invariant vanishes for contact structures with positive \emph{Giroux torsion} \cite{lisstip}.

As the contact invariant is defined in terms of a compatible open book decomposition, we often denote $c(\xi)$ by $c(\Sig, g)$. This class satisfies the following naturality property \cite{osz1}: 

\begin{thm}[Ozsv{\'a}th-Szab{\'o}] 
\label{thm:NaturalityOfContactInvariant}
If $(\Sig, g)$ is an open book decomposition for Y, and $\gamma \subset Y - L$ is a curve supported in a page of the open book (L is the binding), which is not homotopic to the boundary, then $(\Sig,  t_{\gamma}^{-1}g )$ induces an open book decomposition of $Y_{+1}(\gamma)$ (here, $t_{\gamma}$ denotes the right-handed Dehn twist around $\gamma$). And under the map

$$F_{W}: \heeg(-Y) \longrightarrow \heeg(-Y_{+1}(\gamma))$$

\noindent obtained by the two-handle addition (and summing over all $spin^{c}$ structures), we have that 

$$F_{W}(c(\Sig, g)) = \pm c(\Sig, t_{\gamma}^{-1}g).$$

\end{thm}

In general, it is an open question as to whether tightness of the contact structures compatible with open books $(\Sig, g)$ and $(\Sig, h)$ implies tightness of the contact structure compatible with $(\Sig, hg)$. This question is open even when $h$ is a right-handed Dehn twist $t_{\gamma}$ around a curve $\gamma \subset \Sig$ \cite{et2}. Along these lines, however, Theorem \ref{thm:NaturalityOfContactInvariant} implies that if $c(\Sig, g) \neq 0$, then $c(\Sig, t_{\gamma}\circ g) \neq 0$. This behavior is generalized in Theorem \ref{thm:Monoid}, which follows from our main result.
\begin{thm}
\label{thm:Monoid}
Suppose that $g$ and $h$ are diffeomorphisms of $\Sig$ which restrict to the identity on $\partial \Sig$ and that $c(\Sig, g) \neq 0$ and $c(\Sig, h) \neq 0$. Then $c(\Sig, hg) \neq 0$.
\end{thm}

This theorem has the immediate corollary:

\begin{cor}
If the contact structures compatible with $(S,g)$ and $(S,h)$ are strongly fillable, then the contact structure compatible with $(S,hg)$ is tight.
\end{cor}

\noindent Another formulation of Theorem \ref{thm:Monoid} is the statement that for a fixed $\Sig$ the set $$N_S = \{g \in Aut(\Sig, \partial \Sig)\, | \,c(\Sig,g) \neq 0\}$$ is a monoid under composition of diffeomorphisms, where here $Aut(\Sig, \partial \Sig)$ denotes the set of isotopy classes of diffeomorphisms of $\Sig$ that restrict to the identity on $\partial \Sig$. 

The key to this theorem is the observation that the contact invariant satisfies the naturality property below. For any $g$ and $h$ we exhibit a cobordism $X_{g,h}$ with $$\partial X_{g,h} = -Y_{g}-Y_{h} + Y_{hg}$$ where $Y_{\phi}$ is the three-manifold with open book decomposition $(\Sig, \phi)$. This cobordism induces a chain map (multiplication) $$ m: \cheeg(Y_{g}) \x_{\ZZ} \cheeg(Y_{h}) \rightarrow \cheeg(Y_{hg}).$$ If we apply the $\mathrm{Hom}_{\ZZ}(-, \ZZ)$ functor to this expression, we obtain a chain map (comultiplication) $$\mu:\cheeg(-Y_{hg}) \rightarrow \cheeg(-Y_{g}) \x_{\ZZ} \cheeg(-Y_{h}). $$ We show that the contact invariants are natural under the corresponding map $$\tilde{\mu}: \heeg(-Y_{hg}) \rightarrow \heeg(-Y_{g}) \x_{\ZZ} \heeg(-Y_{h})$$ induced on homology. That is 
\begin{thm}
\label{thm:Naturality}
The map $\tilde{\mu}$ takes $c(\Sig, hg) \mapsto c(\Sig, g) \otimes c(\Sig, h)$. 
\end{thm}

Hence, if $c(\Sig, hg) = 0$ then either $c(\Sig, g)=0$ or $c(\Sig, h)=0$, and Theorem \ref{thm:Monoid} follows immediately.

In Section \ref{sec:gen}, we generalize the result of Theorem \ref{thm:Naturality} by examining analogous maps on $\hfp$. We use this generalization to prove the following theorem.

\begin{thm}
\label{thm:U}
If $c^+(\Sig, hg) \in U^d \cdot \hfp(-\Yhg)$, then $c^+(\Sig \#_b \Sig, g \# h) \in U^d \cdot \hfp(-(\Yg \# \Yh)).$ 
\end{thm}

In Theorem \ref{thm:U}, $c^+(\Sig,\phi)$ denotes the image in $\hfp(-Y_{\phi})$ of $c(\Sig, \phi)$ under the natural map $\hf(-Y_{\phi}) \rightarrow \hfp(-Y_{\phi})$. Furthermore, if $(\Sig, g)$ and $(\Sig, h)$ are open books compatible with contact structures $(\Yg,\xi_g)$ and $(\Yh, \xi_h)$, respectively, then $(\Sig \#_b \Sig, g \# h)$ is an open book compatible with the contact structure $(\Yg \# \Yh, \xi_g \# \xi_h),$ where $\#_b$ denotes boundary connected sum.

The following corollaries of Theorem \ref{thm:U} provide obstructions to a contact three-manifold with open book $(\Sig, hg)$ having a compatible open book with planar pages. At the same time, we obtain obstructions to the three-manifold $\Yhg$ being an $L$-space. The reader should compare these corollaries with those found by Ozsv{\'a}th, Szab{\'o}, and Stipsicz in \cite{osz13}. In what follows, $\s(\Sig,\phi)$ denotes the $spin^c$ structure associated to the contact 2-plane field which is compatible with $(\Sig, \phi)$.

\begin{cor}
\label{cor1}
Suppose that $c^+(\Sig \#_b \Sig, g \# h) \neq 0$ and that $c_1(\s(\Sig,g))$ is non-torsion. Then $\Yhg$ is not an $L$-space and the contact structure corresponding to $(\Sig, hg)$ is not compatible with a planar open book.
\end{cor}

\begin{cor}
\label{cor2}
Suppose that $(\Sig, g)$ and $(\Sig, h)$ correspond to Stein fillable contact manifolds $(\Yg, \xi_g)$ and $(\Yh, \xi_h)$ with fillings $(X_g,J_g)$ and $(X_h,J_h).$ Suppose further that $c_1(\s(\xi_g)) =0 = c_1(\s(\xi_h))$ and $c_1(X_g,J_g) \neq 0.$ Then $\Yhg$ is not an $L$-space and the contact structure corresponding to $(\Sig, hg)$ is not compatible with a planar open book.
\end{cor}

\subsection{Acknowledgements}
I would like to thank Ko Honda for bringing to my attention the possibility of Theorem \ref{thm:Monoid}. I also wish to thank John Etnyre, Danny Gillam, Robert Lipshitz, and Shaffiq Welji for very helpful discussions. And, as always, I am indebted to Peter Ozsv{\'a}th for his invaluable comments and suggestions. \\

\section{Heegaard diagrams and the contact class}
\label{sec:Heegaard}
Honda, Kazez, and Mati{\'c} give another interpretation of the Ozsv{\'a}th-Szab{\'o} contact class in \cite{hkm3}. We use their reformulation in our proof of Theorem \ref{thm:Monoid}. Recall that the open book $(\Sig, g)$ is a decomposition of the 3-manifold $Y=\Sig \times [0,1]/ \sim$, where $\sim$ is the relation defined by 
\begin{eqnarray*}
(x,1) \sim(g(x),0), && x \in \Sig\\
(x,t) \sim (x, s), && x\in \partial \Sig, \ t,s \in [0,1]
\end{eqnarray*}

\begin{figure}[!htbp]
\begin{center}
\includegraphics[width=13cm]{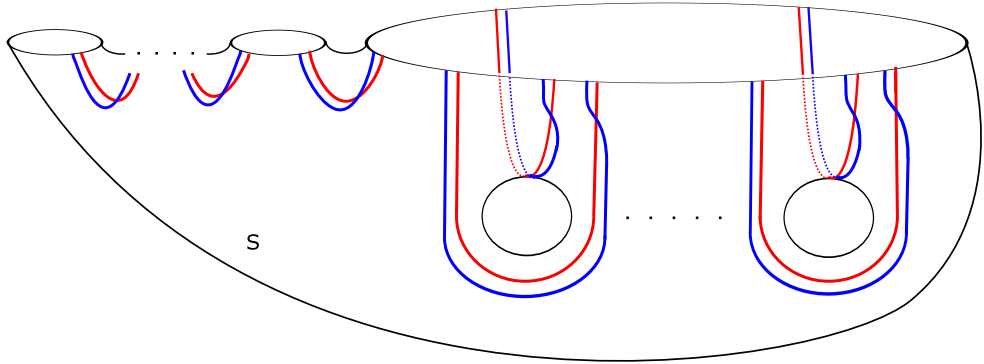}
\caption{\quad A surface $\Sig$ with multiple boundary components and genus greater than one. The $a_i$ arcs are in red and the $b_i$ arcs are in blue.}
\label{fig:ArcExample}
\end{center}
\end{figure}

$Y$ has a Heegaard splitting $Y = H_1 \cup H_2$, where $H_1$ is the handlebody $\Sig \times [0,1/2]$ and $H_2$ is the handlebody $\Sig \times [1/2,1]$. Let $\Sig_t$ denote the page $\Sig \times \{t\}$. If $\Sig$ has $n$ boundary components and genus $k$ then the Heegaard surface in this splitting is the genus $2k+n-1$ surface $\Sigma = \Sig_{1/2} \cup - \Sig_{0}$. To specify a pointed Heegaard diagram for $Y$ it remains to describe the $\alpha$ and $\beta$ curves on $\Sigma$ and the placement of a basepoint $z$. Choose $2k+n-1$ disjoint properly embedded arcs $a_1, \dots, a_{2k+n-1}$ on $\Sigma$ so that $\Sigma \setminus \cup a_i$ is topologically a disk. For each $i$ we obtain $b_i$ by changing the arcs $a_i$ via a small isotopy which moves the endpoints of the $a_i$ along $\partial \Sig$ in the direction given by the orientation of $\partial \Sig$ so that $a_i$ intersects $b_i$ transversely in one point and with positive sign (where $b_i$ inherits its orientation from $a_i$). See Figure \ref{fig:ArcExample} for an illustration of the $a_i$ and $b_i$ arcs on a surface $\Sig$.

\begin{figure}[!htbp]
\begin{center}
\includegraphics[height=12cm]{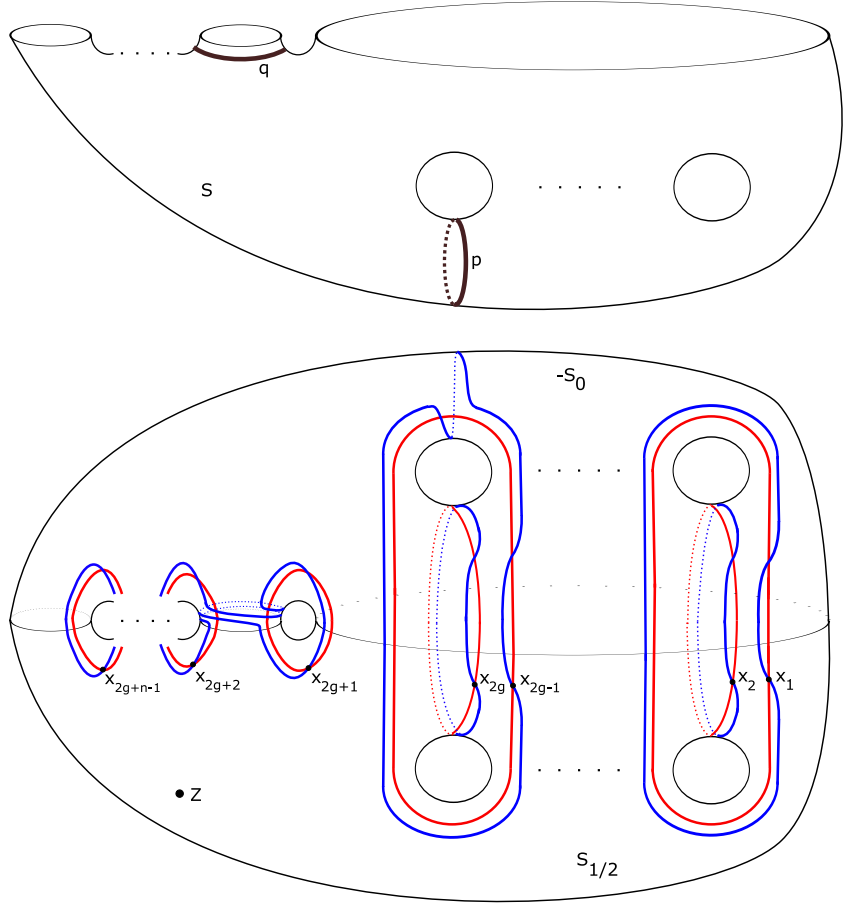}
\caption{\quad The bottom figure is a pointed Heegaard diagram for the open book $(\Sig, g)$, where $\Sig$ is a genus $k$ surface with $n$ boundary components and $g$ is the composition of a left-handed Dehn twist around the curve $p$ with a right-handed Dehn twist around the curve $q$. The $\alpha_i$ are in red and the $\beta_i$ are in blue. We have also labeled the intersection point $\xg = \{x_1,\dots,x_{2k+n-1}\}$.}
\label{fig:HeegExample}
\end{center}
\end{figure}

Now define $$\alpha_i = a_i \times \{1/2\} \cup a_i \times \{0\}$$ $$\beta_i = b_i \times \{1/2\} \cup g(b_i) \times \{0\}.$$ Place the basepoint $z$ in the big region on $\Sig_{1/2}$ (that is, not in one of the thin strip regions). For each $i = 1,\dots,2k+n-1$ let $x_i$ be the intersection point on $\Sig_{1/2}$ between $\alpha_i$ and $\beta_i$. Then $\xg = \{x_1,\dots,x_{2k+n-1}\}$ is an intersection point between $\mathbb{T}_{\alpha}$ and $\mathbb{T}_{\beta}$ in $Sym^{2k+n-1}(\Sigma)$. Moreover, $\xg$ is a cycle in $Hom(\cheeg(Y), \ZZ) = \cheeg(-Y)$ because of the placement of $z$. See Figure \ref{fig:HeegExample} for an illustration of the pointed Heegaard diagram for an open book. 

\begin{thm}[Honda-Kazez-Mati{\'c}] $[\xg] \in \heeg(-Y)$ is the Ozsv{\'a}th-Szab{\'o} contact class $c(\Sig, g)$. 
\end{thm}

\section{Naturality under comultiplication}
\label{sec:Nat}

Given a surface $\Sig$ with genus $k$ and $n$ boundary components, let $a_i$ and $b_i$ be the set of properly embedded arcs described above. We construct another set of disjoint properly embedded arcs $c_i$ from the $b_i$ by changing the arcs $b_i$ via a small isotopy which moves the endpoints of the $b_i$ along $\partial \Sig$ in the direction given by the orientation of $\partial \Sig$. We require that both $a_i$ and $b_i$ intersect $c_i$ transversely in one point and with positive sign (where $c_i$ inherits its orientation from $b_i$). For any two diffeomorphisms $g$ and $h$, we construct three sets of attaching curves on the Heegaard surface $\Sigma = \Sig_{1/2} \cup -\Sig_0$: $$\alpha_i = a_i \times \{1/2\} \cup a_i \times \{0\}$$ $$\beta_i = b_i \times \{1/2\} \cup g(b_i) \times \{0\}$$ $$\gamma_i = c_i \times\{1/2\} \cup h(g(c_i)) \times \{0\}.$$ Once again, we place the basepoint $z$ in the big region of $\Sig_{1/2}$ (outside of the thin strip regions). Then $(\Sigma, \bar{\alpha}, \bar{\beta}, \bar{\gamma}, z)$ is a pointed Heegaard triple-diagram and can be used as in \cite{osz8} to construct a cobordism $X_{\alpha,\beta,\gamma}$ with $$ \partial X_{\alpha,\beta,\gamma} = -Y_{\alpha, \beta} -Y_{\beta,\gamma}+Y_{\alpha,\gamma}$$ where $Y_{\alpha,\beta}$ is the three manifold with Heegaard decomposition $(\Sigma, \bar{\alpha},\bar{\beta})$ (and similarly for $Y_{\beta,\gamma}$ and $Y_{\alpha,\gamma}$). Such a cobordism induces a chain map $$ \cheeg(Y_{\alpha,\beta}) \x_{\ZZ} \cheeg(Y_{\beta,\gamma}) \rightarrow \cheeg(Y_{\alpha,\gamma}).$$ By the description of the Heegaard diagram associated to an open book in section \ref{sec:Heegaard}, it is clear that $$Y_{\alpha,\beta} = Y_{g}$$ $$Y_{\beta,\gamma} = Y_{h}$$ $$Y_{\alpha,\gamma} = Y_{hg}.$$ Thus, we have a chain map 

\begin{eqnarray}
\label{eqn:ChainMap}
m: \cheeg(Y_{g}) \x_{\ZZ} \cheeg(Y_{h}) \rightarrow \cheeg(Y_{hg}).
\end{eqnarray}

If the pointed Heegaard triple-diagram $(\Sigma, \bar{\alpha},\bar{\beta},\bar{\gamma},z)$ is \emph{weakly-admissible} then this map is defined on the generators of $\cheeg(Y_{g}) \x_{\ZZ} \cheeg(Y_{h})$ by $$m(\bar{a}\otimes\bar{b}) = \sum_{\bar{x} \in \mathbb{T}_{\alpha}\cap\mathbb{T}_{\gamma}} \sum_{\{\phi \in \pi_2(\bar{a},\bar{b},\bar{x})\ |\ \mu(\phi)=0,\ n_z(\phi)=0\}} (\# \mathcal{M}(\phi)) \bar{x}.$$ In this sum, $\pi_2(\bar{a},\bar{b},\bar{x})$ is the set of homotopy classes of \emph{Whitney triangles} connecting $\bar{a}$, $\bar{b}$, and $\bar{x}$; $\mu(\phi)$ is the expected dimension of holomorphic representatives of $\phi$; $n_z(\phi)$ is the algebraic intersection number of $\phi$ with the subvariety $\{z\}\times Sym^{2k+n-2}(\Sigma)\subset Sym^{2k+n-1}(\Sigma)$; and $\mathcal{M}(\phi)$ is the moduli space of holomorphic representatives of $\phi$. For more details, see \cite{osz8}.


 As alluded to in the introduction, we apply the $\mathrm{Hom}_{\ZZ}(-,\ZZ)$ functor to the expression in equation \ref{eqn:ChainMap}. If we represent each chain complex diagrammatically by drawing an arrow from $x$ to $y$ whenever $y$ is a term in $\partial x$ or when $y$ is a term in the image of $x$ under the map $m,$ then applying the $\mathrm{Hom}_{\ZZ}(-,\ZZ)$ functor corresponds to reversing the direction of every arrow. Doing so, we obtain a chain map 
 
 \begin{eqnarray*}
\label{eqn:ChainMap2}
\mu: \cheeg(-Y_{hg}) \rightarrow \cheeg(-Y_{g}) \x_{\ZZ} \cheeg(-Y_{h}). 
\end{eqnarray*}

 An element in $ \cheeg(-Y_{hg})$ is a sum of intersection points $\bar{x} \in \mathbb{T}_{\alpha} \cap \mathbb{T}_{\gamma}$. On such an intersection point the chain map $\mu$ is defined by $$\mu(\bar{x}) = \sum_{\bar{a} \in \mathbb{T}_{\alpha}\cap\mathbb{T}_{\beta}, \ \bar{b} \in \mathbb{T}_{\beta}\cap\mathbb{T}_{\gamma}} \sum_{\{\phi \in \pi_2(\bar{a},\bar{b},\bar{x})\ |\ \mu(\phi)=0,\ n_z(\phi)=0\}} (\# \mathcal{M}(\phi)) (\bar{a} \otimes \bar{b})$$ as long as the pointed Heegaard triple-diagram $(\Sigma, \bar{\alpha},\bar{\beta},\bar{\gamma},z)$ is weakly-admissible. To prove Theorem \ref{thm:Naturality}, we show that $$\mu(\overline{x_{hg}}) = \overline{x_ {g}}\otimes\xh.$$ 
 
 We complete the proof in two steps. First we show that $(\Sigma, \bar{\alpha},\bar{\beta},\bar{\gamma},z)$ is weakly-admissible. Then we show that there is only one pair $\{\bar{a}\in \mathbb{T}_{\alpha,\beta}, \ \bar{b} \in \mathbb{T}_{\beta,\gamma}\}$ for which there exists a homotopy class $\phi \in \pi_2(\bar{a},  \bar{b}, \xhg)$ with $n_z(\phi) = 0$ and such that $\phi$ has a holomorphic representative. Moreover, $\bar{a} = \xg$, $\bar{b}=\xh$, and the number of holomorphic representatives of $\phi$ is one.

\subsection{Weak Admissibility} We begin with two definitions from \cite{osz8}.
\begin{definition} For a pointed Heegaard triple-diagram $(\Sigma, \bar{\alpha},\bar{\beta},\bar{\gamma},z)$, let $D_1,\dots,D_m$ be the connected regions of $\Sigma - \cup \alpha_i- \cup \beta_i - \cup \alpha_i$. And let $\psi = \sum_i p_i D_i$ be a formal linear combination of the $D_i$ so that $n_z(\psi) = 0$ and $\partial \psi = \sum_i p_i \partial D_i$ is a linear combination of complete $\alpha$, $\beta$, and $\gamma$ curves. Then $\psi$ is called a \emph{triply-periodic domain}.
\end{definition}

\begin{definition} The pointed Heegaard triple-diagram is said to be \emph{weakly-admissible} if every non-trivial triply-periodic domain $\psi$ has both positive and negative coefficients. 
\end{definition}

\begin{figure}[!htbp]
\begin{center}
\includegraphics[width=13cm]{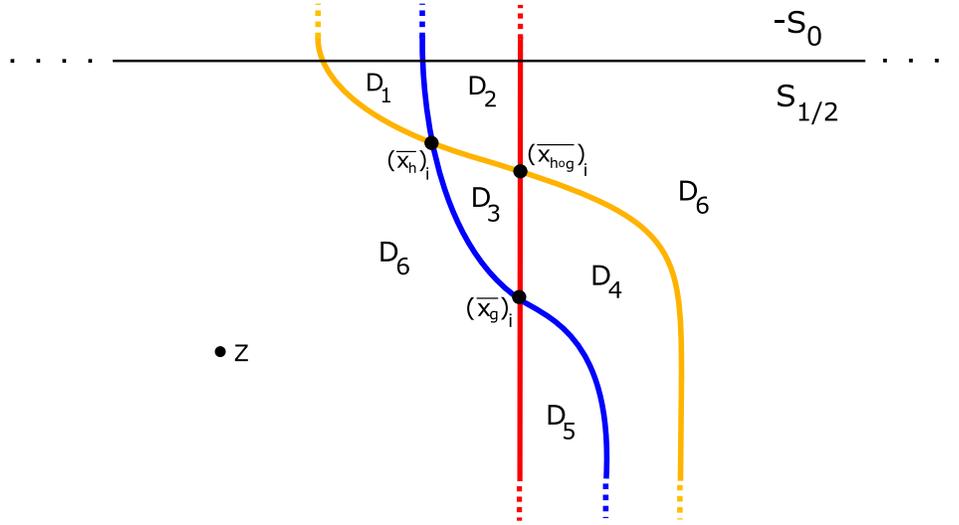}
\caption{\quad The curve $\alpha_i$ is in red, $\beta_i$ in blue, and $\gamma_i$ in yellow. We have labeled regions $D_1,\dots,D_6$ as well as the $ith$ components of the various contact classes $\xg$, $\xh$, and $\xhg$.}
\label{fig:WeakAdm}
\end{center}
\end{figure}

For each $i =1,\dots,2k+n-1$ the curves $\alpha_i$, $\beta_i$, and $\gamma_i$ intersect on $S_{1/2}$ in the arrangement depicted in Figure \ref{fig:WeakAdm}.

If $\psi =\sum_j p_j D_j$ is a triply-periodic domain, then $p_6 = 0$ since $D_6$ contains the basepoint $z$. Since $\partial \psi$ includes some number of complete $\alpha_i$ curves, $$p_2 = p_3-p_4 = -p_5.$$ Therefore, $\psi$ has both positive and negative coefficients unless $p_2=p_5=0$ and $p_3=p_4$. So let's assume the latter. Since $\partial \psi$ includes some number of complete $\beta_i$ curves, $$p_1 = -p_3=0.$$ So, either $\psi$ has both positive and negative coefficients or $p_1 = \dots = p_6 = 0$ and $\partial \psi$ includes no $\alpha_i$, $\beta_i$, or $\gamma_i$ curves. If we carry out this analysis for all $i = 1,\dots, 2k+n-1$ we see that either $\psi$ has both positive and negative coefficients or else it is trivial. Hence $(\Sigma, \bar{\alpha},\bar{\beta},\bar{\gamma},z)$ is weakly-admissible. \qed

\subsection{Completing the proof of Theorem \ref{thm:Naturality}} 
\begin{figure}[!htbp]
\begin{center}
\includegraphics[width=4cm]{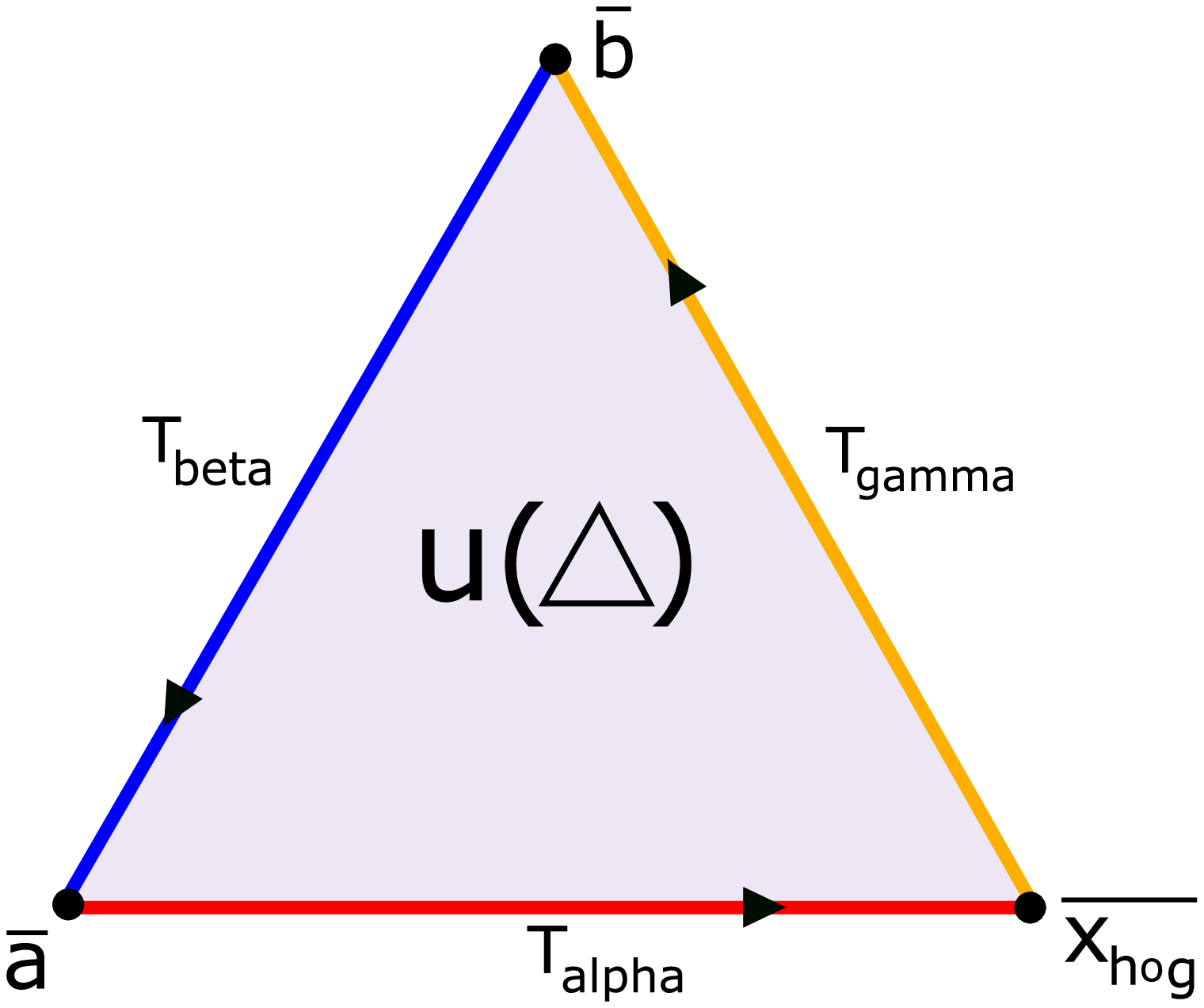}
\caption{\quad}
\label{fig:Whitney}
\end{center}
\end{figure}

 Let $\Delta$ denote the 2-simplex and label its vertices clockwise $v_{\alpha}, v_{\beta}, v_{\gamma}$. Let $e_{\alpha}$ be the edge opposite $v_{\alpha}$ (and similarly for $e_{\beta}$ and $e_{\gamma}$). The boundary of $\Delta$ inherits the standard counterclockwise orientation. Then
\begin{definition} A map $u:\Delta \rightarrow Sym^{2k+n-1}(\Sigma)$ satisfying $u(v_{\gamma})=\bar{a}$, $u(v_{\alpha})=\bar{b}$, and $u(v_{\beta})=\xhg$, and $u(e_{\alpha}) \subset \mathbb{T}_{\alpha}$, $u(e_{\beta}) \subset \mathbb{T}_{\beta}$, and $u(e_{\gamma}) \subset \mathbb{T}_{\gamma}$ is called a \emph{Whitney triangle} between $\bar{a}$, $\bar{b}$, and $\xhg$. This map $u$ is represented schematically in Figure \ref{fig:Whitney}.
\end{definition}

We can represent $\phi \in \pi_2(\bar{a}, \bar{b}, \xhg)$ by a  2-chain $\widehat{\phi} = \sum_j p_j D_j$ whose oriented boundary consists of $\alpha$ arcs from $\bar{a}$ to $\xhg$, $\beta$ arcs from $\bar{b}$ to $\bar{a}$, and $\gamma$ arcs from $\xhg$ to $\bar{b}$. Suppose $n_z(\phi)=0$ and $\phi$ has a holomorphic representative. Then $n_z(\widehat{\phi})=0$ and the $p_j$ are all non-negative. We refer to Figure \ref{fig:WeakAdm} for the local picture near the $ith$ component of the contact classes $\xg$, $\xh$, and $\xhg$. Write $$\widehat{\phi} = p_1D_1+\dots +p_6D_6+\sum_{j>6}p_jD_j.$$ 

Now we can analyze the possibilities for $p_1,\dots,p_6$ given the boundary constraints on $\widehat{\phi}$. $\xhgi$ must be a corner of the region defined by $\widehat{\phi}$; moreover this corner is such that we enter $\xhgi$ along an arc of $\alpha_i$ and we leave along an arc of $\gamma_i$. Therefore, $p_6+p_3=p_2+p_4+1$. If $\xhi$ is not a corner, then $p_3+p_1=p_2+p_6$. Note that $p_6 = 0$ since $n_z(\widehat{\phi})$ = 0. Thus, these two equations become 

 \begin{eqnarray*}
p_3&=&p_2+p_4+1 \\
p_3+p_1&=&p_2.
\end{eqnarray*}

Subtracting the second equation from the first, we have $$-p_1=p_4+1$$ which implies that either $p_1$ or $p_4$ is negative, which cannot happen since $\phi$ has a holomorphic representative. Therefore, $\xhi$ \emph{is} a corner. The same type of analysis shows that $\xgi$ is a corner. 

Since $\xhi$ is a corner, either $p_1+p_3+1 = p_2$ or $p_1+p_3 = p_2+1$. Substituting $p_3 = p_2+p_4+1$ into both expressions, we have the two possibilities $p_1+p_2+p_4+2 = p_2$ or $p_1+p_2+p_4=p_2$. We can rule out the first possibility as it implies that either $p_1$ or $p_4$ is negative. And the second possibility holds only if $p_1 = p_4 = 0$. So, to summarize what we know so far: $p_1=0$, $p_3=p_2+1$, $p_4=0$, and $p_6=0.$ 

Since $\xgi$ is a corner, then either $p_5+p_3=p_4+p_6+1$ or $p_5+p_3+1=p_4+p_6$. Substituting what we know of $p_3$, $p_4$, and $p_6$ into these two expressions, we obtain the two possibilities $p_5+p_2+1 = 1$ or $p_5+p_2+2 = 0$. We can rule out the second possibility as it implies that either $p_5$ or $p_2$ is negative. And the first possibility holds only if $p_5 = p_2 = 0$. Thus, we have determined that the only possibility for the values $p_1,\dots,p_6$ are:  $$p_1=p_2=p_4=p_5=p_6=0$$ $$p_3=1.$$ 

Because the same analysis works for every $i = 1,\dots,2k+n-1$ and because every component of $\partial \widehat{\phi}$ must contain some $\xhgi$ we can conclude that $\widehat{\phi}$ is the linear combination which is the sum of precisely one of these small triangular regions ($D_3$ in figure \ref{fig:WeakAdm}) for each $i$. Therefore, any holomorphic triangle $\phi$ between $\bar{a}$, $\bar{b}$, and $\xhg$ with $n_z(\phi)=0$ is, in fact, a triangle between $\xg$, $\xh$, and $\xhg$, and can be expressed as a product of these small triangles in our Heegaard diagram. Moreover, since each of these disjoint triangular regions is topologically a disk, and we have specified the image of three boundary points, $\#\mathcal{M}(\phi)=1$ by the Riemann Mapping Theorem. Hence, $$\mu(\overline{x_{hg}}) = \overline{x_ {g}}\otimes\xh.$$ Therefore, $$\tilde{\mu}(c(\Sig, hg)) = c(\Sig,g)\otimes c(\Sig,h)$$ and the proof of Theorem \ref{thm:Naturality} is complete. As mentioned in the Introduction, Theorem \ref{thm:Monoid} follows immediately. \qed \\ 

\section{A generalization of Theorem \ref{thm:Naturality} }
\label{sec:gen}
\subsection{$\cfl$ and connected sums}
\label{ssec:CFL}
For a $spin^c$ structure $\s$ on $Y$ and a pointed Heegaard diagram $(\Sigma, \al, \be, z)$ for $Y$ which is \emph{strongly $\s$-admissible}, recall that we can define a chain complex $\cfl(Y,\s)$ which is finitely generated as a $\ZZ[U]$ module \cite{osz8}. The generators of $\cfl(Y,\s)$ are pairs of the form $[\bar{x},i]$, where $\bar{x} \in \mathbb{T}_{\alpha} \cap \mathbb{T}_{\beta}$, $\s_z(\bar{x}) = \s$, $i \in \mathbb{Z}^{\leq 0}$, and $U$ acts by $U \cdot [\bar{x},i] = [\bar{x},i-1]$. The differential on $\cfl(Y, \s)$ is given by $$\partial [\bar{x},i] = \sum_{\bar{y} \in \mathbb{T}_{\alpha}\cap\mathbb{T}_{\gamma}} \sum_{\{\phi \in \pi_2(\bar{x},\bar{y})\ |\ \mu(\phi)=1\}} \#( \frac{\mathcal{M}(\phi)}{\mathbb{R}}) \cdot [\bar{y}, i-n_z(\phi)].$$

We can identify $\cf(Y,\s)$ with $\cfl(Y,\s)/U\cdot \cfl(Y, \s)$, so there is a natural quotient map $$\pi: \cfl(Y,\s) \rightarrow \cf(Y,\s).$$ If $A$ is a $\ZZ[U]$ module, let $A^{\lor}$ denote $\mathrm{Hom}_{\ZZ[U]}(A, \ZZ[U,U^{-1}]/\ZZ[U])$. Then observe that we can identify $\cfl(Y,\s)^{\lor}$ with $\cfp(-Y,\s)$, and $\cf(Y,\s)^{\lor}$ with $\cf(-Y,\s)$ as complexes over $\ZZ[U]$ and $\ZZ$, respectively. \footnote{We have already been identifying $\cf(-Y)$ with $\mathrm{Hom}_{\ZZ}(\cf(Y), \ZZ)$, but the latter is isomorphic as a $\ZZ$ module to $\cf(Y)^{\lor}$ since we are thinking of $\cf(Y)$ as a $\ZZ[U]$ module where the action of $U$ is multiplication by 0.} In fact, applying the $\mathrm{Hom}_{\ZZ[U]}(-, \ZZ[U,U^{-1}]/\ZZ[U])$ functor to the expression above, we obtain the natural inclusion map $$\pi^{\lor}: \cf(-Y,\s) \rightarrow \cfp(-Y,\s)$$ which sends $$\bar{x} \mapsto [\bar{x},0]$$ for $\bar{x} \in \mathbb{T}_{\alpha} \cap \mathbb{T}_{\beta}$. 

In \cite{osz14} the authors construct a homotopy equivalence $$f^{\leq 0}: \cfl(Y_1,\s_1)\x_{\ZZ[U]} \cfl(Y_2,\s_2) \rightarrow \cfl(Y_1 \# Y_2,\s_1 \# \s_2)$$ as follows. Let $(\Sigma_i,\al_i,\be_i, z_i)$ be a strongly $\s_i$-admissible pointed Heegaard diagram for $Y_i$, for $i = 1, \, 2$. Consider the pointed Heegaard triple-diagram $(\Sigma_1 \# \Sigma_2, \al_1 \cup \al_2', \be_1 \cup \al_2,\be_1' \cup \al_2, z)$, where the connected sum of $\Sigma_1$ with $\Sigma_2$ occurs at the $z_i$, $z$ is a point in the connected sum region, and $\al_2'$ and $\be_1'$ are exact Hamiltonian translates of $\al_2$ and $\be_1$ so that this new diagram is admissible. Suppose that the genera of $\Sigma_1$ and $\Sigma_2$ are $g_1$ and $g_2$ and let $\Theta_1 \in \mathbb{T}_{\beta_1} \cap \mathbb{T}_{\beta_1'}$ and $\Theta_2 \in \mathbb{T}_{\alpha_2} \cap \mathbb{T}_{\alpha_2'}$ be the top graded intersection points in $\cf(\#^{g_1}(S^1 \times S^2), \s_0)$ and $\cf(\#^{g_2}(S^1 \times S^2), \s_0)$. Then the maps $[\bar{x},i] \mapsto [\bar{x} \times \Theta_2, i]$ and $[\bar{y},j] \mapsto [\Theta_1 \times \bar{y},j]$ define chain maps $$\Phi_1:\cfl(Y_1, \s_1) \rightarrow \cfl(Y_1\# \#^{g_2} (S^1 \times S^2),\s_1 \# \s_0),$$ and $$\Phi_2:\cfl(Y_2, \s_2) \rightarrow \cfl( \#^{g_2} (S^1 \times S^2)\#Y_2,\s_0 \# \s_2).$$ The pointed Heegaard triple-diagram above can be used to define a map 

\begin{eqnarray*}
\Gamma: \cfl(Y_1\# \#^{g_2} (S^1 \times S^2),\s_1 \# \s_0) \x_{\ZZ[U]}\cfl( \#^{g_2} (S^1 \times S^2)\#Y_2,\s_0 \# \s_2)&& \\
\longrightarrow \cfl(Y_1 \# Y_2, \s_1 \# \s_2).&&
\end{eqnarray*}

The map $f^{\leq 0}$ is then defined to be $$f^{\leq0} = \Gamma \circ (\Phi_1 \x \Phi_2).$$ The maps in this composition are all $U$-equivariant by construction; hence, so is $f^{\leq 0}.$ In \cite{osz14}, Oszv{\'a}th and Szab{\'o} show that $f^{\leq 0}$ is a homotopy equivalence. We define a homotopy equivalence $$f: \cf(Y_1,\s_1)\x_{\ZZ} \cf(Y_2,\s_2) \rightarrow \cf(Y_1 \# Y_2, \s_1 \# \s_2)$$ in exactly the same way. \footnote{In \cite{osz14}, the authors define a homotopy equivalence between these two chain complexes in a slightly different and more direct way. However, the map $f$ defined here is better suited for our purposes.}

\subsection{$\hfp$ and the contact invariant}
Let $\phi$ denote the Whitney triangle between $\xg$, $\xh$, and $\xhg$ found in Section \ref{sec:Nat}, and let $\s(\phi)$ denote the $spin^c$ structure on $\Xgh$ corresponding to $\phi$. Moreover, let $\sg$, $\sh$, and $\shg$ denote the induced $spin^c$ structures on $\Yg$, $\Yh$, and $\Yhg$. Then there is a chain map $$m_{\s(\phi)}: \cf(Y_{g},\sg) \x_{\ZZ} \cf(Y_{h},\sh) \rightarrow \cf(Y_{hg},\shg)$$ which is a refinement of the map $m$ defined in Section \ref{sec:Nat}. The difference between the two is that $m_{\s(\phi)}$ counts only those Whitney triangles which correspond to the $spin^c$ structure $\s(\phi)$. Applying the $\mathrm{Hom}_{\ZZ[U]}(-, \ZZ[U,U^{-1}]/\ZZ[U])$ functor as before and taking homology, it is clear that the induced map $$(m_{\s(\phi)}^{\lor})_*: \hf(-\Yhg,\shg) \rightarrow \hf(-\Yg,\sg) \x_{\ZZ} \hf(-\Yh,\sh)$$ still takes $c(\Sig, hg) \mapsto c(\Sig, g) \x c(\Sig, h).$ Mirroring the notation in Section \ref{sec:Nat}, we denote this map by $\tilde{\mu}_{\s(\phi)}$.

The pointed Heegaard triple-diagram $(\Sigma, \al, \be, \ga, z)$ from Section \ref{sec:Nat} also gives a $U$-equivariant chain map $$m^{\leq 0}_{\s(\phi)}: \cfl(Y_{g},\sg) \x_{\ZZ[U]} \cfl(Y_{h},\sh) \rightarrow \cfl(Y_{hg},\shg).\footnote{For this map to be well-defined, we need $(\Sigma, \al, \be, \ga, z)$ to be strongly $\s(\phi)$-admissible. We can assume, however, that this is the case since the Heegaard diagram $(\Sigma, \al, \be, \ga, z)$ can be made strongly $\s(\phi)$-admissible by winding the $\al$, $\be$, and $\ga$ circles  around curves contained strictly in the $-\Sig_0$ portion of the surface $\Sigma$ \cite{osz8}.}$$ Let $$\pi_{g \x h}: \cfl(Y_{g},\sg) \x_{\ZZ[U]} \cfl(Y_{h},\sh) \rightarrow \cf(Y_{g},\sg) \x_{\ZZ} \cf(Y_{h},\sh)$$ and $$\pi_{hg}:\cfl(Y_{hg},\shg)\rightarrow \cf(Y_{hg},\shg)$$ denote the quotient maps discussed in Subsection \ref{ssec:CFL}. Then the following diagram commutes.

\begin{displaymath}
\xymatrix{
\cfl(Y_{g},\sg) \x_{\ZZ[U]} \cfl(Y_{h},\sh)    \ar[dd]^{\pi_{g \x h}}  \ar[rrr]^{ \hspace{14mm}     m^{\leq 0 }_{\s(\phi)}}&&& \cfl(Y_{hg},\shg) \ar[dd]^{\pi_{hg}}\\
&&&\\
\cf(Y_{g},\sg) \x_{\ZZ} \cf(Y_{h},\sh)  \ar[rrr]^{   \hspace{14mm}m_{\s(\phi)}}         &&&\cf(Y_{hg},\shg)      }
\end{displaymath}

After applying the $\mathrm{Hom}_{\ZZ[U]}(-, \ZZ[U,U^{-1}]/\ZZ[U])$ functor and taking homology, we obtain the commutative diagram

\begin{displaymath}
\xymatrix{
  H_*((\cfl(Y_{g},\sg) \x_{\ZZ[U]} \cfl(Y_{h},\sh))^{\lor})  &&& \ar[lll]^{\hspace{14mm}\tilde{\mu}^+_{\s(\phi)}} \hfp(-Y_{hg},\shg)  \\
&&&\\
     \hf(-Y_{g},\sg) \x \hf(-Y_{h},\sh) \ar[uu]^{(\pi_{g \x h}^{\lor})_*}   &&& \ar[lll]^{\hspace{14mm}\tilde{\mu}_{\s(\phi)}} \ar[uu]^{(\pi_{hg}^{\lor})_*}  \hf(-Y_{hg},\shg)      }
\end{displaymath}
where $\tilde{\mu}^+_{\s(\phi)} = ((m^{\leq 0}_{\s(\phi)})^{\lor})_*$.
Recall that for a contact three-manifold $(Y,\xi)$, the class $c^+(\xi) \in \hfp(-Y)$ is defined to be the image of $c(\xi)$ under the natural map $\hf(-Y) \rightarrow \hfp(-Y)$ \cite{osz13}. As was mentioned in Subsection \ref{ssec:CFL}, $(\pi_{hg}^{\lor})_*$ is this natural map, and therefore $$(\pi_{hg}^{\lor})_*( c(\Sig, hg))= c^+(\Sig, hg).$$ Let $K_{g \x h}$ denote $\tilde{\mu}^+_{\s(\phi)} (c^+(\Sig, hg)).$ Then by the commutativity of this diagram, $$K_{g \x h}=\tilde{\mu}^+_{\s(\phi)} (c^+(\Sig, hg))=(\pi_{g \x h}^{\lor})_*(c(\Sig, g) \x c(\Sig, h)).$$

Returning to our discussion of connected sums, let $$\pi_{g \#h}:\cfl(\Yg \# \Yh, \sg \# \sh)\rightarrow \cf(\Yg \# \Yh, \sg \# \sg)$$ denote the natural quotient map, and let $f^{\leq 0}$ and $f$ be the homotopy equivalences described in Subsection \ref{ssec:CFL}. Then the diagram below commutes.

\begin{displaymath}
\xymatrix{
\cfl(Y_{g},\sg) \x_{\ZZ[U]} \cfl(Y_{h},\sh)    \ar[dd]^{\pi_{g \x h}}  \ar[rrr]^{\hspace{8mm}f^{\leq 0}}&&& \cfl(\Yg \# \Yh, \sg \# \sh) \ar[dd]^{\pi_{g \#h}}\\
&&&\\
\cf(Y_{g},\sg) \x_{\ZZ} \cf(Y_{h},\sh)  \ar[rrr]^{\hspace{8mm}f}         &&&\cf(\Yg \# \Yh, \sg \# \sh)    }
\end{displaymath}

Again, after applying the $\mathrm{Hom}_{\ZZ[U]}(-, \ZZ[U,U^{-1}]/\ZZ[U])$ functor and taking homology, we obtain the commutative diagram

\begin{displaymath}
\xymatrix{
  H_*((\cfl(Y_{g},\sg) \x_{\ZZ[U]} \cfl(Y_{h},\sh))^{\lor})  &&& \ar[lll]^{\hspace{10mm}\tilde{k}^+} \hfp(-(\Yg \# \Yh),\sg \# \sh)  \\
&&&\\
     \hf(-Y_{g},\sg) \x_{\ZZ} \hf(-Y_{h},\sh) \ar[uu]^{(\pi_{g \x h}^{\lor})_*}   &&& \ar[lll]^{\hspace{10mm}\tilde{k}} \ar[uu]^{(\pi_{g \# h}^{\lor})_*}  \hf(-(\Yg \# \Yh), \sg \# \sh)      }
\end{displaymath} where $\tilde{k}^+ = ((f^{\leq 0})^{\lor})_*$ and $\tilde{k} = (f^{\lor})_*$. At this point, we are ready to generalize Theorem \ref{thm:Naturality}.

\begin{thm}
\label{thm:Gen}
The map $(\tilde{k}^+)^{-1} \circ \tilde{\mu}^+_{\s(\phi)}: \hfp(-Y_{hg},\shg) \rightarrow \hfp(-(\Yg \# \Yh),\sg \# \sh)$ is $U$-equivariant and takes $c^+(\Sig, hg) \mapsto c^+(\Sig \#_b \Sig, g \# h)$.
\end{thm}

Observe that Theorem \ref{thm:U} follows immediately as a corollary. Recall from the introduction that the contact structure compatible with the open book $(\Sig \#_b \Sig, g \# h)$ is the connected sum of the contact structures compatible with $(\Sig,g)$ and $(\Sig,h)$. 

\begin{proof}[Proof of Theorem \ref{thm:Gen}]
The maps $\tilde{\mu}^+_{\s(\phi)}$ and $(\tilde{k}^+)^{-1}$ are certainly $U$-equivariant, so their composition is as well. Moreover, $\tilde{k} = (f^{\lor})_*$ takes $c(\Sig \#_b \Sig, g \# h) \mapsto c(\Sig, g) \x c(\Sig, h)$. This follows from precisely the same sort of argument as was used in Section \ref{sec:Nat} to show that $\tilde{\mu} = (m^{\lor})_*$ takes $c(\Sig, hg) \mapsto c(\Sig, g) \x c(\Sig, h)$. In the pointed Heegaard triple-diagram $(\Sigma_1 \# \Sigma_2, \al_1 \cup \al_2', \be_1 \cup \al_2,\be_1' \cup \al_2, z)$ used to construct the map $\Gamma$, the only holomorphic Whitney triangle $\phi \in \pi_2(\bar{a},\bar{b},\xg \times \xg)$ with $n_z(\phi) = 0$ is a product of small triangles connecting $\xg \times \Theta_2$, $\Theta_1 \times \xh$, and $\xg \times \xh$. Therefore, the map $f^{\lor}$ takes $\xg \times \xh \mapsto \xg \x \xh$ on the level of chains, and $(f^{\lor})_*$ takes $c(\Sig \#_b \Sig, g \# h) \mapsto c(\Sig, g) \x c(\Sig, h)$.

Hence, $$((\pi_{g \x h}^{\lor})_* \circ \tilde{k}) (c(\Sig \#_b \Sig, g \# h) ) = (\pi_{g \x h}^{\lor})_*(c(\Sig, g) \x c(\Sig, h))=K_{g \x h}.$$ Therefore, by commutativity, $$K_{g \x h}=(\tilde{k}^+ \circ (\pi_{g \# h}^{\lor})_*)(c(\Sig \#_b \Sig, g \# h)) = \tilde{k}^+(c^+(\Sig \#_b \Sig, g \# h)).$$ Thus, $(\tilde{k}^+)^{-1}(K_{g \x h}) = c^+(\Sig \#_b \Sig, g \# h).$ But this implies that  $$((\tilde{k}^+)^{-1} \circ \tilde{\mu}^+_{\s(\phi)} )(c^+(\Sig, hg)) = (\tilde{k}^+)^{-1}(K_{g \x h}) = c^+(\Sig \#_b \Sig, g \# h).$$
\end{proof}

\section{$L$-spaces and planar open books}
\label{sec:planar}
Etnyre recently showed that while every overtwisted contact structure has a compatible open book with planar pages, there are fillable contact structures which do not \cite{et}. More recently, Ozsv{\'a}th, Szab{\'o}, and Stipsicz found a Heegaard Floer homology obstruction to a contact structure having a compatible open book with planar pages, and were able to reproduce some of Etnyre's results \cite{osz13}. Their main result is the following:

\begin{thm}[Ozsv{\'a}th-Szab{\'o}-Stipsicz]
\label{thm:OSS}
If the contact three-manifold $(Y,\xi)$ has a compatible open book with planar pages then $c^+(\xi) \in U^d \cdot \hfp(-Y)$ for all $d \in \mathbb{N}$.
\end{thm}

They prove the following corollaries which are based upon this principle.

\begin{cor}[Ozsv{\'a}th-Szab{\'o}-Stipsicz]
\label{coro1}
Suppose that $c^+(\xi) \neq 0$ and $c_1(\s(\xi))$ is non-torsion. Then it cannot be the case that $c^+(\xi) \in U^d \cdot \hfp(-Y)$ for all $d \in \mathbb{N}$. In particular, $(Y,\xi)$ is not supported by a planar open book.
\end{cor}

\begin{cor}[Ozsv{\'a}th-Szab{\'o}-Stipsicz]
\label{coro2}
Suppose that the contact three-manifold $(Y,\xi)$ has a Stein filling $(X,J)$ with $c_1(\s(\xi))=0$ and $c_1(X,J) \neq 0$. Then it cannot be the case that $c^+(\xi) \in U^d \cdot \hfp(-Y)$ for all $d \in \mathbb{N}$. In particular, $(Y,\xi)$ is not supported by a planar open book.
\end{cor}

Our Corollaries \ref{cor1} and \ref{cor2} now follow from Theorems \ref{thm:U}, \ref{thm:Gen}, and the above corollaries of Ozsv{\'a}th, Szab{\'o}, and Stipsicz.

\begin{proof}[Proof of Corollary \ref{cor1}]
If $c_1(\s(\Sig, g))$ is non-torsion, then so is $$c_1(\s(\Sig \#_b \Sig, g \# h)) = c_1(\s(\Sig,g)) \oplus c_1(\s(\Sig, h)).$$ If, in addition, $c^+(\Sig \#_b \Sig, g \# h) \neq 0$, then Corollary \ref{coro1} implies that it cannot be the case that $c^+(\Sig \#_b \Sig, g \# h) \in  U^d \cdot \hfp(-(\Yg \# \Yh))$ for all $d \in \mathbb{N}$. Thus, Theorem \ref{thm:U} demands that it cannot be the case that $c^+(\Sig, hg)\in  U^d \cdot \hfp(-\Yhg)$ for all $d \in \mathbb{N}$. It follows that $\Yhg$ cannot be an $L$-space. Furthermore, it follows from Theorem \ref{thm:OSS} of Oszv{\'a}th, Szab{\'o}, and Stipsicz that the contact structure supported by $(\Sig, hg)$ is not compatible with a planar open book.

\end{proof}

\begin{proof}[Proof of Corollary \ref{cor2}]
If $(\Sig, g)$ and $(\Sig, h)$ correspond to contact manifolds $(Y_g,\xi_g)$ and $(Y_h, \xi_h)$ with Stein fillings $(X_g, J_g)$ and $(X_h, J_h)$, then $(\Sig \#_b \Sig, g \# h)$ corresponds to the contact manifold $(Y_g \# Y_h, \xi_g \# \xi_h)$ with Stein filling $(X_g \#_b X_h, J_g \#_b J_h)$. If $c_1(\s(\xi_g)) = 0$ and $c_1(\s(\xi_h)) = 0$, then $$c_1(\s(\xi_g \# \xi_h)) =  c_1(\s(\xi_g)) \oplus c_1(\s(\xi_h)) = 0.$$ Moreover, if $c_1(X_g, J_g) \neq 0$, then $$c_1(X_g \#_b X_h, J_g \#_b J_h) = c_1(X_g, J_g) \oplus c_1(X_h, J_h) \neq 0.$$ Then, by Corollary \ref{coro2}, it cannot be the case that $c^+(\Sig \#_b \Sig, g \# h) \in  U^d \cdot \hfp(-(\Yg \# \Yh))$ for all $d \in \mathbb{N}$. And, just as before, Theorem \ref{thm:U} then implies that it cannot be the case that $c^+(\Sig, hg)\in  U^d \cdot \hfp(-\Yhg)$ for all $d \in \mathbb{N}$. It follows that $\Yhg$ cannot be an $L$-space and that the contact structure supported by $(\Sig, hg)$ is not compatible with a planar open book.
\end{proof}

It is not clear whether we can replace the condition that $c^+(\Sig \#_b \Sig, g \# h) \neq 0$ in the formulation of Corollary \ref{cor1} by the condition that $c^+(\Sig,g) \neq 0$ and $c^+(\Sig, h) \neq 0$. For the contact invariant $c$ defined in $\hf(-Y)$, $c(\Sig, g) \neq 0$ and $c(\Sig, h) \neq 0$ implies that $c(\Sig \#_b \Sig, g \# h) \neq 0$. It is not immediately obvious that the same is true in general for the contact invariant $c^+$. There are, however, \emph{special} cases in which the same holds for $c^+$. For instance, if $(\Sig, g)$ and $(\Sig, h)$ support strongly fillable contact structures, then so does $(\Sig \#_b \Sig, g \# h)$, and hence, $c^+(\Sig \#_b \Sig, g \# h) \neq 0$. More useful perhaps, is the following.

\begin{claim} 
\label{claim:NZ}
If $c^+(\Sig, g) \neq 0$ and $(\Sig, h)$ supports a Stein fillable contact structure, then $c^+(\Sig \#_b \Sig, g \# h)\neq 0$.
\end{claim}

\begin{proof}[Proof of Claim \ref{claim:NZ}]
If $(\Sig, h)$ supports a Stein fillable contact structure, then after a number of positive stabilizations $(\Sig, h)$ is equivalent to an open book $(\Sig', \phi)$, where $\phi$ is the composition of right-handed Dehn twists around curves in $\Sig'$. Hence, $(\Sig \#_b \Sig, g \# h)$ is equivalent via the same number of positive stabilizations to the open book $(\Sig \#_b \Sig', g \# \phi).$ Therefore, by the naturality of the contact invariant $c^+$ under composition with left-handed Dehn twists, if $c^+(\Sig \#_b \Sig', g \# id) \neq 0$, then $c^+(\Sig \#_b \Sig', g \# \phi) \neq 0$. \footnote{Theorem \ref{thm:Naturality} was stated only for the contact invariant $c$, but the same holds for $c^+$ via a commutative diagram chase.} Let $m = 2k+r$, where $k$ is the genus of $\Sig'$ and $r$ is the number of boundary components of $\Sig'$. Then $(\Sig \#_b \Sig', g \# id)$ is an open book for the manifold $\Yg \# \#^{m}(S^1 \times S^2)$. There is an isomorphism $$\hfp(-(\Yg \# \#^{m}(S^1 \times S^2)), \s_g \# \s_0) \rightarrow \hfp(-\Yg,\s_g) \land^* H^1(\#^{m}(S^1 \times S^2))$$ \cite[Proposition 6.4]{osz14}. It is clear from the construction that this isomorphism takes $$c^+(\Sig \#_b \Sig', g \# id) \mapsto c^+(\Sig,g) \land \Theta$$ where $\Theta$ is the lowest graded element of $ \land^* H^1(\#^{m}(S^1 \times S^2)).$ Therefore, $c^+(\Sig, g) \neq 0$ if and only if $c^+(\Sig \#_b \Sig', g \# id) \neq 0$. Putting these facts together, we get the claim.
\end{proof}

It remains to be seen what can be shown using these techniques. It would be very interesting, as Etnyre mentions in \cite{et}, to find a non-fillable contact structure which is not supported by an open book with planar pages. To this end it is enough to find, by Corollary \ref{cor1}, a Stein fillable open book $(\Sig, g)$ with $c_1(\s(\Sig, g))$ non-torsion and an open book $(\Sig, h)$ with $c^+(\Sig,h) \neq 0$ such that $(\Sig, hg)$ is non-fillable. 
 
\bibliographystyle{hplain.bst}
\bibliography{References}

\end{document}